%% file: paper.tex
\documentclass[12pt]{article}

\usepackage{amsfonts, amssymb, amsmath, amsthm, epsf}
\makeatletter

\@addtoreset{figure}{section}

\theoremstyle{plain}
\newtheorem{theorem}{Theorem}
\@addtoreset{theorem}{section}

\newtheorem{lemma}{Lemma}
\@addtoreset{lemma}{section}

\@addtoreset{corollary}{section}

\newtheorem{condition}{Condition}
\@addtoreset{condition}{section}

\theoremstyle{definition}

\newtheorem{definition}{Definition}
\@addtoreset{definition}{section}

\newtheorem{example}{Example}
\@addtoreset{example}{section}

\newtheorem{remark}{Remark}
\@addtoreset{remark}{section}

\@addtoreset{equation}{section}

\makeatother

\newcommand{\loc}{{\rm loc}}

\newcommand{\supp}{{\rm supp\,}}
\newcommand{\Dom}{{\rm Dom\,}}
\renewcommand{\ker}{{\rm ker\,}}
\renewcommand{\dim}{{\rm dim\,}}
\newcommand{\codim}{{\rm codim\,}}
\newcommand{\Gr}{{\rm Gr\,}}

\newcommand{\Orb}{{\rm Orb}}

\newcommand{\dist}{{\rm dist}}
\newcommand{\Span}{{\rm span}}

\title{Generalized Solutions\\ of Nonlocal Elliptic Problems}
\author{Pavel Gurevich\thanks{This work was supported by
Russian Foundation for Basic Research (grant No.~02-01-00312),
Russian Ministry for Education (grant No.~E02-1.0-131), and INTAS
(grant~YSF~2002-008).}}

\date{}

\begin{document}

\maketitle

\begin{abstract}
An elliptic equation of order $2m$ with general nonlocal
boundary-value conditions, in a plane bounded domain $G$ with
piecewise smooth boundary, is considered. Generalized solutions
belonging to the Sobolev space $W_2^m(G)$ are studied. The
Fredholm property of the unbounded operator corresponding to the
elliptic equation, acting on $L_2(G)$, and defined for functions
from the space $W_2^m(G)$ that satisfy homogeneous nonlocal
conditions is proved.
\end{abstract}

\input introd.tex
\input sect1.tex
\input sect2.tex

\input bibl.tex
\end{document}

%% file: introd.tex
\section*{Introduction}
In the one-dimensional case, nonlocal problems were studied by
A.~Sommerfeld~\cite{Sommerfeld}, J.~D.~Tamarkin~\cite{Tamarkin},
M.~Picone~\cite{Picone}. T.~Carleman~\cite{Carleman} considered
the problem of finding a function harmonic on a two-dimensional
bounded domain and subjected to a nonlocal condition connecting
the values of this function at different points of the boundary.
A.~V.~Bitsadze and A.~A.~Smarskii~\cite{BitzSam} suggested another
setting of a nonlocal problem arising in plasma theory: to find a
function harmonic on a two-dimensional bounded domain and
satisfying nonlocal conditions on shifts of the boundary that can
take points of the boundary inside the domain. Different
generalizations of the above nonlocal problems were investigated
by many authors~\cite{Vishik,Browder,ZhEid,RSh,Kishk,GM,Gusch}.

It turns out that the most difficult situation occurs if the
support of nonlocal terms intersects the boundary. In that case,
solutions of nonlocal problems can have power-law singularities
near some points even if the boundary and the right-hand sides are
infinitely smooth~\cite{SkMs86, SkRJMP}. For this reason, such
problems are naturally studied in weighted spaces (introduced by
V.~A.~Kondrat'ev for boundary-value problems in nonsmooth
domains~\cite{KondrTMMO67}). The most complete theory of nonlocal
problems in weighted spaces is developed by
A.~L.~Skubachevskii~\cite{SkMs86,SkDu90,SkDu91,SkJMAA,SkBook}.

Note that the investigation of nonlocal problems is motivated both
by significant theoretical progress in that direction and
important applications arising in biophysics, theory of diffusion
processes~\cite{Feller}, plasma theory~\cite{Sam}, and so on.

\smallskip

In the present paper, we study {\it generalized solutions} of an
elliptic equation of order $2m$ in a two-dimensional bounded
domain $G$, satisfying nonlocal boundary-value conditions that are
set on parts $\Gamma_j$ of the boundary $\partial
G=\bigcup_j\overline{\Gamma_j}$. By generalized solutions, we mean
functions from the Sobolev space $W^m(G)=W_2^m(G)$. We prove that
an unbounded operator acting on $L_2(G)$ and corresponding to the
above nonlocal problem has the Fredholm property.

Note that solutions of nonlocal problems can be sought on the
space of ``smooth'' functions, namely, on the Sobolev space
$W^{2m}(G)$ (see~\cite{GurRJMP03,GurRJMP04}) or on weighted spaces
$H_a^{2m}(G)$, where $$
 \|u\|_{H_a^k(G)}=
 \bigg(\sum_{|\alpha|\le k}\int_G \rho^{2(a-k+|\alpha|)}
 |D^\alpha u|^2\bigg)^{1/2},
$$ $k\ge 0$ is an integer, $a\in\mathbb R,$
$\rho=\rho(y)=\dist(y,\mathcal K)$, and $\mathcal
K=\bigcup_j\overline{\Gamma_j}\setminus\Gamma_j$ is the set formed
by finitely many points of conjugation of nonlocal conditions
(see~\cite{SkMs86,SkDu91}). In both cases, a {\it bounded}
operator corresponds to the nonlocal problem. Whether or not this
operator has the Fredholm property depends on spectral properties
of some auxiliary problems with a parameter. In turn, these
spectral properties are affected by the values of the coefficients
in nonlocal conditions and by a geometrical structure of the
support of nonlocal terms and the boundary near the set $\mathcal
K$. However, if we consider generalized solutions (i.e., functions
from $W^m(G)$), then the corresponding {\it unbounded} operator
turns out to have the Fredholm property irrespective of the above
factors.

Earlier the Fredholm property of an unbounded nonlocal operator on
$L_2(G)$ was studied either for the case in which nonlocal
conditions were set on shifts of the boundary~\cite{SkBook} or in
the case of a nonlocal perturbation of the Dirichlet problem for a
second-order elliptic equation~\cite{GM,Gusch}. Elliptic equations
of order $2m$ with general nonlocal conditions are investigated
for the first time.

%% file: sect1.tex
\section{Setting of Nonlocal Problems in Bounded Domains}\label{sectStatement}

\subsection{Setting of Problem}\label{subsectStatement}
Let $G\subset{\mathbb R}^2$ be a bounded domain with boundary
$\partial G$. Consider a set ${\mathcal K}\subset\partial G$
consisting of finitely many points. Let $\partial
G\setminus{\mathcal K}=\bigcup\limits_{i=1}^{N}\Gamma_i$, where
$\Gamma_i$ are open (in the topology of $\partial G$)
$C^\infty$-curves. We assume that, in a neighborhood of each point
$g\in{\mathcal K}$, the domain $G$ is a plane angle.

Denote by ${\bf P}(y, D_y)={\bf P}(y,D_{y_1},D_{y_2})$ and
$B_{i\mu s}(y, D_y)=B_{i\mu s}(y,D_{y_1},D_{y_2})$ differential
operators of order $2m$ and $m_{i\mu}$ ($m_{i\mu}\le m-1$),
respectively, with complex-valued $C^\infty$ coefficients, and let
${\bf P}^0(y, D_y)$ and $B^0_{i\mu s}(y, D_y)$ denote their
principal homogeneous parts ($i=1, \dots, N;$ $\mu=1, \dots, m;$
$s=0, \dots, S_i$). Here $D_y=(D_{y_1},D_{y_2})$,
$D_{y_j}=-i\partial/\partial y_j$.

Now we formulate conditions on the operators ${\bf P}(y, D_y)$ and
$B_{i\mu 0}(y, D_y)$ (these operators will correspond to a
``local'' elliptic problem). We assume that the operator ${\bf
P}(y, D_y)$ is {\it properly elliptic} on $\overline{G}$; in
particular, the following estimate holds for all $\theta\in\mathbb
R^2$ and $y\in\overline{G}$:
\begin{equation}\label{eqPEllipinG}
 A^{-1}|\theta|^{2m}\le |{\bf P}^0(y, \theta)|\le
 A|\theta|^{2m},\quad A>0.
\end{equation}
Further, let $y\in\overline{\Gamma_i}$. One may assume with no
loss of generality that the curve $\overline{\Gamma_i}$ is defined
by the equation $y_2=0$ near the point $y$. We suppose that the
{\it system $\{B_{i\mu 0}(y, D_y)\}_{\mu=1}^m$ satisfies the
Lopatinsky condition with respect to the operator ${\bf P}(y,
D_y)$ for all $i=1, \dots, N$}. In other words, let the polynomial
$$
B_{i\mu 0}'(y, \tau)\equiv\sum\limits_{\nu=1}^m
b_{i\mu\nu}(y)\tau^{\nu-1}\equiv B_{i\mu 0}^0(y,1,\tau)\ \big({\rm
mod\,} \mathbf M^+(y,\tau)\big)
$$
be the residue of dividing $B_{i\mu 0}^0(y, 1,\tau)$ by $\mathbf
M^+(y,\tau)$, where
$$
\mathbf M^+(y,\tau)=\prod\limits_{\nu=1}^m(\tau-\tau^+_\nu(y)),
$$
while $ \tau^+_1(y),\dots, \tau^+_m(y)$ are the roots of the
polynomial ${\bf P}^0(y, 1, \tau)$ with positive imaginary parts
(note that ${\bf P}^0(y, 1, \tau)$, $B_{i\mu 0}^0(y, 1,\tau)$, and
$\mathbf M^+(y,\tau)$ are considered as polynomials in $\tau$). In
this case, the validity of the Lopatinsky condition means that
$$
 d_i(y)=\det\|b_{i\mu\nu}(y)\|_{\mu,\nu=1}^m\ne0.
$$
Since each of the curves $\overline{\Gamma_i}$, $i=1,\dots,N$, is
a compact, it follows that
\begin{equation}\label{eqCover}
 D=\min\limits_{i=1,\dots,N}\,\inf\limits_{y\in\overline{\Gamma_i}}|d_i(y)|>0.
\end{equation}
We emphasize that the operators $B_{i\mu 0}(y, D_y)$ are not
necessarily normal on $\overline{\Gamma_i}$.

\medskip

For an integer $k\ge0$, denote by $W^k(G)=W_2^k(G)$ the Sobolev
space with the norm
$$
\|u\|_{W^k(G)}=\left(\sum\limits_{|\alpha|\le k}\int\limits_G
|D^\alpha u|^2\,dy\right)^{1/2}
$$
(we set $W^0(G)=L_2(G)$ for $k=0$). For an integer $k\ge1$, we
introduce the space $W^{k-1/2}(\Gamma)$ of traces on a smooth
curve $\Gamma\subset\overline{ G}$ with the norm
\begin{equation}\label{eqTraceNormW}
\|\psi\|_{W^{k-1/2}(\Gamma)}=\inf\|u\|_{W^k(G)}\quad (u\in W^k(G):
u|_\Gamma=\psi).
\end{equation}

\smallskip

Denote $\mathbf B_{i\mu}^0 u=B_{i\mu 0}(y, D_y)u(y)|_{\Gamma_i}$.
As we have mentioned above, the operators ${\bf P}(y, D_y)$ and
$\mathbf B_{i\mu}^0$ will \textit{correspond to a ``local''
boundary-value problem}.

\smallskip

Now we define operators corresponding to nonlocal conditions near
the set $\mathcal K$. Let $\Omega_{is}$ ($i=1, \dots, N;$ $s=1,
\dots, S_i$) be $C^\infty$-diffeomorphisms taking some
neighborhood ${\mathcal O}_i$ of the curve
$\overline{\Gamma_i\cap\mathcal O_{{\varepsilon}}(\mathcal K)}$
onto the set $\Omega_{is}({\mathcal O}_i)$ in such a way that
$$
\Omega_{is}(\Gamma_i\cap\mathcal O_{{\varepsilon}}(\mathcal
K))\subset G,
$$
\begin{equation}\label{eqOmega}
\Omega_{is}(g)\in\mathcal K\qquad\text{for}\qquad
g\in\overline{\Gamma_i}\cap\mathcal K.
\end{equation}
Here ${\varepsilon}>0$, $\mathcal O_{\varepsilon}(\mathcal
K)=\{y\in \mathbb R^2: \dist(y, \mathcal K)<\varepsilon\}$ is the
$\varepsilon$-neighborhood of the set $\mathcal K$. Thus, under
the transformations $\Omega_{is}$, the curves
$\Gamma_i\cap\mathcal O_{{\varepsilon}}(\mathcal K)$ are mapped
strictly inside the domain $G$, whereas the set of end points
$\overline{\Gamma_i}\cap\mathcal K$ is mapped to itself.

Let us specify the structure of the transformations $\Omega_{is}$
near the set $\mathcal K$. Denote by the symbol $\Omega_{is}^{+1}$
the transformation $\Omega_{is}:{\mathcal
O}_i\to\Omega_{is}({\mathcal O}_i)$ and by $\Omega_{is}^{-1}$ the
transformation $\Omega_{is}^{-1}:\Omega_{is}({\mathcal
O}_i)\to{\mathcal O}_i$ inverse to $\Omega_{is}.$ The set of all
points
$$
\Omega_{i_qs_q}^{\pm1}(\dots\Omega_{i_1s_1}^{\pm1}(g))\in{\mathcal
K}\qquad (1\le s_j\le S_{i_j},\ j=1, \dots, q),
$$
i.e., the set of all points that can be obtained by consecutively
applying the transformations $\Omega_{i_js_j}^{+1}$ or
$\Omega_{i_js_j}^{-1}$ (taking the points of ${\mathcal K}$ to
${\mathcal K}$) to the point $g\in\mathcal K$ is called an {\it
orbit} of $g$ and is denoted by $\Orb(g)$.

Clearly, for any $g, g'\in{\mathcal K}$ either $\Orb(g)=\Orb(g')$
or $\Orb(g)\cap\Orb(g')=\varnothing$. In what follows, we suppose
that the set $\mathcal K$ consists of a unique orbit. (All results
can be directly generalized for the case in which $\mathcal K$
consists of finitely many mutually disjoint orbits.) The set
(orbit) $\mathcal K$ consists of $N$ points, which we denote by
$g_j$, $j=1, \dots, N$.

Take a small number $\varepsilon$ (see Remark~\ref{remSmallEps}
below) such that there exist neighborhoods $\mathcal
O_{\varepsilon_1}(g_j)$ of the points $g_j\in\mathcal K$
satisfying the following conditions:
\begin{enumerate}
\item
$ \mathcal O_{\varepsilon_1}(g_j)\supset\mathcal
O_{\varepsilon}(g_j) $;
\item in the neighborhood $\mathcal O_{\varepsilon_1}(g_j)$, the boundary $\partial G$
is a plane angle;
\item
$\overline{\mathcal O_{\varepsilon_1}(g_j)}\cap\overline{\mathcal
O_{\varepsilon_1}(g_k)}=\varnothing$ for any $g_j,g_k\in\mathcal
K$, $k\ne j$;
\item if $g_j\in\overline{\Gamma_i}$ and
$\Omega_{is}(g_j)=g_k,$ then ${\mathcal
O}_{\varepsilon}(g_j)\subset\mathcal
 O_i$ and
 $\Omega_{is}\big({\mathcal
O}_{\varepsilon}(g_j)\big)\subset{\mathcal
O}_{\varepsilon_1}(g_k).$
\end{enumerate}

For each point $g_j\in\overline{\Gamma_i}\cap\mathcal K$, we fix a
transformation $y\mapsto y'(g_j)$ of the argument; this
transformation is the composition of the shift by the vector
$-\overrightarrow{Og_j}$ and a rotation by some angle such that
the set ${\mathcal O}_{\varepsilon_1}(g_j)$ is mapped onto the
neighborhood ${\mathcal O}_{\varepsilon_1}(0)$ of the origin,
while the sets
$$
G\cap{\mathcal
O}_{\varepsilon_1}(g_j)\qquad\text{and}\qquad\Gamma_i\cap{\mathcal
O}_{\varepsilon_1}(g_j)
$$
are mapped onto the intersection of a plane angle
$$
K_j=\{y\in{\mathbb R}^2:\ r>0,\ |\omega|<\omega_j\}
$$
with the neighborhood ${\mathcal O}_{\varepsilon_1}(0)$ and the
intersection of  the side
$$
\gamma_{j\sigma}=\{y\in\mathbb R^2:\
\omega=(-1)^\sigma \omega_j\}
$$
($\sigma=1$ or $\sigma=2$) of the angle $K_j$ with the
neighborhood ${\mathcal O}_{\varepsilon_1}(0)$, respectively. Here
$(\omega,r)$ are the polar coordinates of the point $y$ and
$0<\omega_j<\pi$.

\begin{condition}\label{condK1}
The above change of variables $y\mapsto y'(g_j)$ for
$y\in{\mathcal O}_{\varepsilon}(g_j)$,
$g_j\in\overline{\Gamma_i}\cap\mathcal K$, reduces the
transformation $\Omega_{is}(y)$ {\rm(}$i=1,\dots, N,$ $s=1,\dots,
S_i${\rm)} to the composition of a rotation and a homothety in the
new variables $y'$.
\end{condition}

\begin{remark}\label{remK1}
In particular, Condition~\ref{condK1} combined with the assumption
$\Omega_{is}(\Gamma_i)\subset G$ means that, if
$g\in\Omega_{is}(\overline{\Gamma_i}\cap\mathcal
K)\cap\overline{\Gamma_j}\cap{\mathcal K}\ne\varnothing$, then the
curves $\Omega_{is}(\overline{\Gamma_i})$ and
$\overline{\Gamma_j}$ are not tangent to each other at the point
$g$.
\end{remark}

Consider a number $\varepsilon_0$,
$0<\varepsilon_0\le\varepsilon$, satisfying the following
condition: if $g_j\in\overline{\Gamma_i}$ and
$\Omega_{is}(g_j)=g_k,$ then ${\mathcal
O}_{\varepsilon_0}(g_k)\subset \Omega_{is}\big({\mathcal
O}_{\varepsilon}(g_j)\big)$. Introduce a function $\zeta\in
C^\infty(\mathbb R^2)$ such that
\begin{equation}\label{eqZeta}
 \zeta(y)=1\quad\text{for}\quad y\in\mathcal O_{\varepsilon_0/2}(\mathcal K),\qquad
 \supp\zeta\subset\mathcal O_{\varepsilon_0}(\mathcal K).
\end{equation}

Now we define nonlocal operators $\mathbf B_{i\mu}^1$ by the
formula
$$
 \mathbf B_{i\mu}^1u=\sum\limits_{s=1}^{S_i}
   \big(B_{i\mu s}(y,
   D_y)(\zeta u)\big)\big(\Omega_{is}(y)\big)\qquad\text{for}\qquad
   y\in\Gamma_i\cap\mathcal O_{\varepsilon}(\mathcal K),
$$
$$
 \mathbf B_{i\mu}^1u=0\qquad\text{for}\qquad
y\in\Gamma_i\setminus(\Gamma_i\cap\mathcal
O_{\varepsilon}(\mathcal K)),
$$
where $\big(B_{i\mu s}(y,
D_y)u\big)\big(\Omega_{is}(y)\big)=B_{i\mu s}(x,
D_{x})u(x)|_{x=\Omega_{is}(y)}$. Since $\mathbf B_{i\mu}^1u=0$
whenever $\supp u\subset\overline{ G}\setminus\overline{\mathcal
O_{\varepsilon_0}(\mathcal K)}$, we say that the operators
$\mathbf B_{i\mu}^1$ \textit{correspond to nonlocal terms
supported near the set} $\mathcal K$.

\smallskip

For any $\rho>0$, we denote $G_\rho=\{y\in G: \dist(y,
\partial G)>\rho\}$. Consider operators $\mathbf
B_{i\mu}^2$ satisfying the following condition
(cf.~\cite{SkMs86,SkJMAA,GurRJMP03}).
\begin{condition}\label{condSeparK23}
There exist numbers $\varkappa_1>\varkappa_2>0$ and $\rho>0$ such
that the inequalities
\begin{equation}\label{eqSeparK23'}
  \|\mathbf B^2_{i\mu}u\|_{W^{2m-m_{i\mu}-1/2}(\Gamma_i)}\le c_1
  \|u\|_{W^{2m}(G\setminus\overline{\mathcal O_{\varkappa_1}(\mathcal
  K)})},
\end{equation}
\begin{equation}\label{eqSeparK23''}
  \|\mathbf B^2_{i\mu}u\|_{W^{2m-m_{i\mu}-1/2}
   (\Gamma_i\setminus\overline{\mathcal O_{\varkappa_2}(\mathcal K)})}\le
  c_2 \|u\|_{W^{2m}(G_\rho)}
\end{equation}
hold for any
$$
 u\in W^{2m}(G\setminus\overline{\mathcal
O_{\varkappa_1}(\mathcal
  K)})\cap W^{2m}(G_\rho),
$$
where $i=1, \dots, N$, $\mu=1, \dots, m$, and $c_1,c_2>0$.
\end{condition}

It follows from~\eqref{eqSeparK23'} that $\mathbf B_{i\mu}^2u=0$
whenever $\supp u\subset \mathcal O_{\varkappa_1}(\mathcal K)$.
For this reason, we say that the operators $\mathbf B_{i\mu}^2$
\textit{correspond to nonlocal terms supported outside the set}
$\mathcal K$.

\smallskip

We will suppose throughout that Conditions~\ref{condK1}
and~\ref{condSeparK23} hold.

\smallskip

We study the following nonlocal elliptic problem:
\begin{align}
 {\bf P}(y, D_y)u=f_0(y) \quad &(y\in G),\label{eqPinG}\\
     \mathbf B_{i\mu}u\equiv\mathbf B_{i\mu}^0 u+\mathbf B_{i\mu}^1 u+\mathbf B_{i\mu}^2 u=
   0\quad
    &(y\in \Gamma_i;\ i=1, \dots, N;\ \mu=1, \dots, m),\label{eqBinG}
\end{align}
where $f_0\in L_2(G)$. Introduce the space $W_B^m(G)$ consisting
of functions $u\in W^m(G)$ that satisfy homogeneous nonlocal
conditions~(\ref{eqBinG}): $\mathbf B_{i\mu}u=0$.

Consider the unbounded operator $\mathbf P: \Dom(\mathbf P)\subset
L_2(G)\to L_2(G)$ given by
$$
 \mathbf P u=\mathbf P(y, D_y)u,\qquad u\in \Dom(\mathbf P)=\{u\in
 W_B^{m}(G): \mathbf P(y, D_y) u\in L_2(G)\}.
$$

\begin{definition}\label{defGenSol1}
A function $u$ is called a {\it generalized solution} of
problem~(\ref{eqPinG}), (\ref{eqBinG}) with right-hand side
$f_0\in L_2(G)$ if $u\in \Dom(\mathbf P)$ and $
 \mathbf P u=f_0.
$
\end{definition}

One can give another (equivalent) definition for a generalized
solution. To do so, we write the operator ${\bf P}(y, D_y)$ in the
divergent form,
$$
  {\bf P}(y, D_y)=\sum\limits_{0\le|\xi|,|\beta|\le m}
    D^{\beta}p_{\xi\beta}(y)D^{\xi},
$$
where $p_{\xi\beta}$ are infinitely differentiable functions.

For any set $X\in \mathbb R^2$ having a nonempty interior, denote
by $C_0^\infty(X)$ the set of functions infinitely differentiable
on $\overline{ X}$ and compactly supported on $X$.

\begin{definition}\label{defGenSol2}
A function $u$ is called a {\it generalized solution} of
problem~(\ref{eqPinG}), (\ref{eqBinG}) with right-hand side
$f_0\in L_2(G)$ if $u\in W_B^{m}(G)$ and the integral identity
$$
 \sum\limits_{0\le|\xi|,\,|\beta|\le m}\,\int\limits_G
    p_{\xi\beta}(y)D^{\xi}u\overline{D^{\beta}v}\,dy=\int\limits_G
    f_0\overline{v}\,dy
$$
holds for any $v\in C_0^\infty(G)$.
\end{definition}

\begin{remark}
Generalized solutions a priori belong to the space $W^m(G)$,
whereas Condition~\ref{condSeparK23} is formulated for functions
belonging to the space $W^{2m}$ inside the domain and near the
smooth part of the boundary. Such a formulation can be justified
by the fact that any generalized solution belongs to $W^{2m}$
outside an arbitrarily small neighborhood of the set $\mathcal K$
(see Lemma~\ref{lSmoothOutsideK} below).
\end{remark}

\begin{remark}\label{remSmallEps}
We have supposed above that the number $\varepsilon$ is small
(whereas $\varkappa_1,\varkappa_2,\rho$ can be arbitrary). Let us
show that this leads to no loss of generality. Let us have a
number $\hat\varepsilon$, $0<\hat\varepsilon<\varepsilon$. Take a
number $\hat\varepsilon_0$,
$0<\hat\varepsilon_0\le\hat\varepsilon$, satisfying the following
condition: if $g_j\in\overline{\Gamma_i}$ and
$\Omega_{is}(g_j)=g_k,$ then ${\mathcal
O}_{\hat\varepsilon_0}(g_k)\subset \Omega_{is}\big({\mathcal
O}_{\hat\varepsilon}(g_j)\big)$. Consider a function $\hat\zeta\in
C^\infty(\mathbb R^2)$ such that $ \hat\zeta(y)=1$ for
$y\in\mathcal O_{\hat\varepsilon_0/2}(\mathcal K)$ and
$\supp\hat\zeta\subset\mathcal O_{\hat\varepsilon_0}(\mathcal K)$.
Introduce the operators $\mathbf B_{i\mu}^1$ as follows:
$$
\mathbf B_{i\mu}^1u=\sum\limits_{s=1}^{S_i}
   \big(B_{i\mu s}(y,
   D_y)(\zeta u)\big)\big(\Omega_{is}(y)\big)\qquad\text{for}\qquad
   y\in\Gamma_i\cap\mathcal O_{\varepsilon}(\mathcal K),
$$
$$
  \mathbf
B_{i\mu}^1u=0\qquad\text{for}\qquad
y\in\Gamma_i\setminus(\Gamma_i\cap\mathcal
O_{\varepsilon}(\mathcal K)).
$$
Clearly,
$$
\mathbf B_{i\mu}^0+\mathbf B_{i\mu}^1+\mathbf B_{i\mu}^2= \mathbf
B_{i\mu}^0+\hat{\mathbf B}_{i\mu}^1+\hat{\mathbf B}_{i\mu}^2,
$$
where $\hat{\mathbf B}_{i\mu}^2=\mathbf B_{i\mu}^1-\hat{\mathbf
B}_{i\mu}^1+\mathbf B_{i\mu}^2$. Since $\mathbf
B_{i\mu}^1u-\hat{\mathbf B}_{i\mu}^1u=0$ near the set $\mathcal
K$, it follows that the operator $\mathbf B_{i\mu}^1-\hat{\mathbf
B}_{i\mu}^1$ satisfies Condition~\ref{condSeparK23} for some
suitable $\varkappa_1,\varkappa_2,\rho$
(see~\cite[\S~1]{GurRJMP03} for more details). Thus, we can always
choose $\varepsilon$ to be as small as necessary (possibly at the
expense of a modification of the operator $\mathbf B_{i\mu}^2$ and
the values of $\varkappa_1,\varkappa_2,\rho$).
\end{remark}

\subsection{Example of Nonlocal Problem}
One can consider the following example as a model one.
\begin{example}\label{exGeneralProblem}
Let ${\bf P}(y, D_y)$ and $B_{i\mu s}(y, D_y)$ be the same
operators as above. Let $\Omega_{is}$ ($i=1, \dots, N;$
$s=1,\dots, S_i$) be $C^\infty$-diffeomorphisms taking some
neighborhood ${\mathcal O}_i$ of the (\textit{whole}) curve
$\Gamma_i$ to the set $\Omega_{is}({\mathcal O}_i)$ in such a way
that $\Omega_{is}(\Gamma_i)\subset G$. Consider the following
nonlocal problem:
\begin{equation}\label{eqPinGEx}
 {\bf P}(y, D_y) u=f_0(y) \quad (y\in G),
\end{equation}
\begin{equation}\label{eqBinGEx}
 \begin{aligned}
 B_{i\mu 0}(y, D_y)u(y)|_{\Gamma_i}+
 \sum\limits_{s=1}^{S_i}
   \big(B_{i\mu s}(y, D_y)u\big)\big(\Omega_{is}(y)\big)\big|_{\Gamma_i}=
   0\\
    (y\in \Gamma_i;\ i=1, \dots, N;\ \mu=1, \dots, m).
 \end{aligned}
\end{equation}
We emphasize that \textit{a priori the transformations
$\Omega_{is}$ are not supposed to satisfy
condition~\eqref{eqOmega}$;$ however, we further represent the
nonlocal operators as the sum of the operators $\mathbf
B_{i\mu}^0$, $\mathbf B_{i\mu}^1$, and $\mathbf B_{i\mu}^2$, and
the transformations occurring in the definition of the operators
$\mathbf B_{i\mu}^1$ will satisfy condition~\eqref{eqOmega}.} To
obtain this representation, we take a small ${\varepsilon}$ such
that, for any point $g\in\mathcal K$, the set $\overline{\mathcal
O_{{\varepsilon}}(g)}$ intersects the curve
$\overline{\Omega_{is}(\Gamma_i)}$ only if
$g\in\overline{\Omega_{is}(\Gamma_i)}$. If
$g\in\overline{\Gamma_i}\cap\mathcal K$ and
$\Omega_{is}(g)\in\mathcal K$, then we assume that the
transformation $\Omega_{is}(y)$ satisfies Condition~\ref{condK1}
for $y\in\mathcal O_\varepsilon(g)$.

\begin{remark}
By Remark~\ref{remK1}, Condition~\ref{condK1} is a restriction on
the geometrical structure of the support of nonlocal terms near
the set $\mathcal K$. However, if
$\Omega_{is}(\overline{\Gamma_i}\cap\mathcal K)\subset\partial
G\setminus{\mathcal K}$, then we impose no restrictions on the
geometrical structure of the curve
$\Omega_{is}(\overline{\Gamma_i})$ near the boundary $\partial G$
(cf.~\cite{SkMs86, SkDu91}).
\end{remark}

Let $\zeta\in C^\infty(\mathbb R^2)$ be a function satisfying
relations~\eqref{eqZeta}. Introduce the operators $$\mathbf
B_{i\mu}^0 u=B_{i\mu 0}(y, D_y)u(y)|_{\Gamma_i},$$ $$ \mathbf
B_{i\mu}^1u=\sum\limits_{s=1}^{S_i} \big(B_{i\mu s}(y, D_y)(\zeta
u)\big)\big(\Omega_{is}(y)\big)\big|_{\Gamma_i},$$ $$\mathbf
B_{i\mu}^2u=\sum\limits_{s=1}^{S_i}
   \big(B_{i\mu s}(y,
   D_y)((1-\zeta)
   u)\big)\big(\Omega_{is}(y)\big)\big|_{\Gamma_i}
$$ (see figures~\ref{figB1} and~\ref{figB2}).
\begin{figure}[p]
{ \hfill\epsfbox{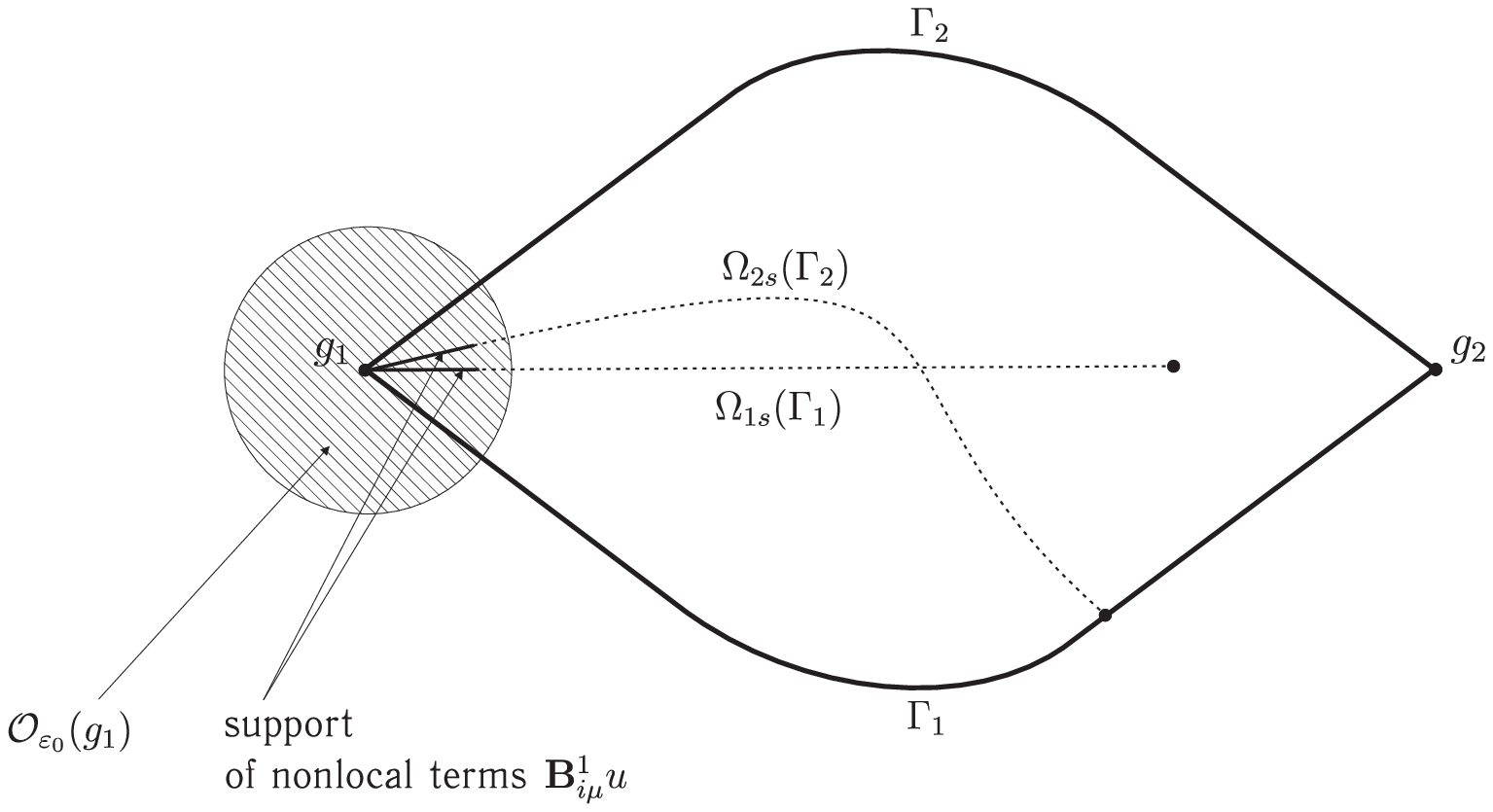}\hfill\ } \caption{To
problem~\eqref{eqPinGEx}, \eqref{eqBinGEx}}
   \label{figB1}
\end{figure}
\begin{figure}[p]
{ \hfill\epsfbox{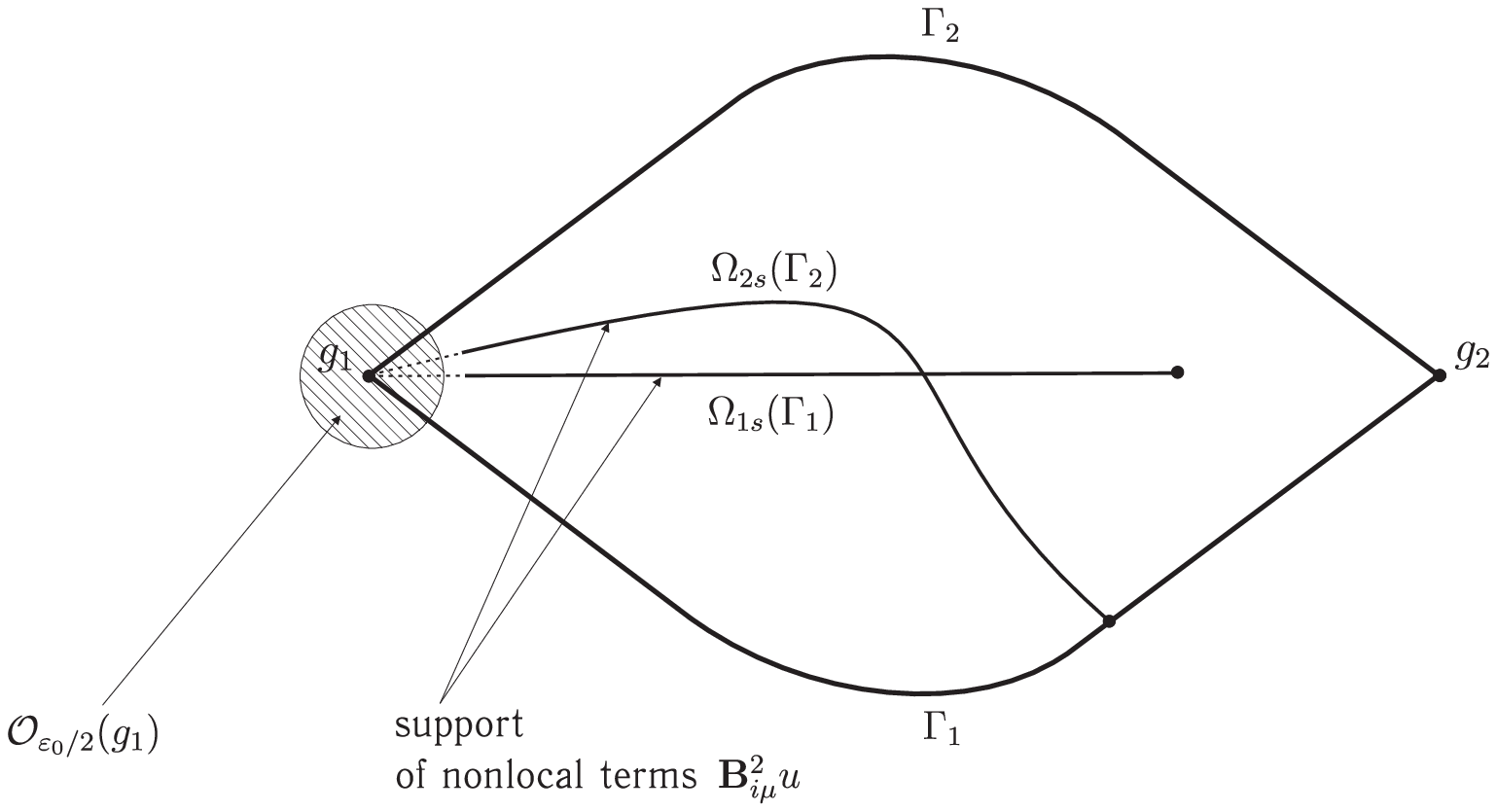}\hfill\ } \caption{To
problem~\eqref{eqPinGEx}, \eqref{eqBinGEx}}
   \label{figB2}
\end{figure}
Since the support of the function $\zeta$ is concentrated near the
set $\mathcal K$, one may assume that the transformations
$\Omega_{is}$ occurring in the definition of the operator $\mathbf
B_{i\mu}^1$ are defined on some neighborhood of the set $\mathcal
K$ and satisfy condition~\eqref{eqOmega}. Moreover, it follows
from~\cite[Sec.~1.2]{GurRJMP03} that the operator $\mathbf
B_{i\mu}^2$ satisfies Condition~\ref{condSeparK23}. Therefore,
problem~(\ref{eqPinGEx}), (\ref{eqBinGEx}) can be represented in
the form~(\ref{eqPinG}), (\ref{eqBinG}).
\end{example}

\subsection{Nonlocal Problems near the Set $\mathcal
K$}\label{subsectStatementNearK}

When studying problem~(\ref{eqPinG}), (\ref{eqBinG}), one must pay
special attention to the behavior of solutions near the set
${\mathcal K}$ of conjugation points. Now we consider the
corresponding model problems.

Denote by $u_j(y)$ the function $u(y)$ for $y\in{\mathcal
O}_{\varepsilon_1}(g_j)$. If $g_j\in\overline{\Gamma_i},$
$y\in{\mathcal O}_{\varepsilon}(g_j),$ and
$\Omega_{is}(y)\in{\mathcal O}_{\varepsilon_1}(g_k),$ then we
denote the function $u(\Omega_{is}(y))$ by $u_k(\Omega_{is}(y))$.
In this notation, nonlocal problem~(\ref{eqPinG}), (\ref{eqBinG})
acquires the following form in the $\varepsilon$-neighborhood of
the set (orbit) $\mathcal K$:
\begin{gather*}
 {\bf P}(y, D_y) u_j=f_0(y) \quad (y\in\mathcal O_\varepsilon(g_j)\cap
 G),\\
\begin{aligned}
B_{i\mu 0}(y, D_y)u_j(y)|_{\mathcal
O_\varepsilon(g_j)\cap\Gamma_i}+ \sum\limits_{s=1}^{S_i}
\big(B_{i\mu s}(y,D_y)(\zeta
u_k)\big)\big(\Omega_{is}(y)\big)\big|_{\mathcal
O_\varepsilon(g_j)\cap\Gamma_i}
=f_{i\mu}(y) \\
\big(y\in \mathcal O_\varepsilon(g_j)\cap\Gamma_i;\ i\in\{1\le
i\le N: g_j\in\overline{\Gamma_i}\};\ j=1, \dots, N;\ \mu=1,
\dots, m\big),
\end{aligned}
\end{gather*}
where $f_{i\mu}=-\mathbf B_{i\mu}^2u$.

Let $y\mapsto y'(g_j)$ be the change of variables described in
Sec.~\ref{subsectStatement}. Denote
$K_j^\varepsilon=K_j\cap\mathcal O_\varepsilon(0)$ and
$\gamma_{j\sigma}^\varepsilon=\gamma_{j\sigma}\cap\mathcal
O_\varepsilon(0)$. Introduce the functions
$$
U_j(y')=u_j(y(y')),\qquad f_j(y')=f_0(y(y')),\qquad y'\in
K_j^\varepsilon,\qquad
$$
$$
f_{j\sigma\mu}(y')=f_{i\mu}(y(y')),\qquad
y'\in\gamma_{j\sigma}^\varepsilon,
$$
where $\sigma=1$ $(\sigma=2)$ if, under the transformation
$y\mapsto y'(g_j)$, the curve $\Gamma_i$ is mapped to the side
$\gamma_{j1}$ ($\gamma_{j2}$) of the angle $K_j$. Denote $y'$ by
$y$ again. Then, by virtue of Condition~\ref{condK1},
problem~(\ref{eqPinG}), (\ref{eqBinG}) acquires the form
\begin{gather}
  {\bf P}_{j}(y, D_y)U_j=f_{j}(y) \quad (y\in
  K_{j}^\varepsilon),\label{eqPinK}\\
 \sum\limits_{k,s}
       (B_{j\sigma\mu ks}(y, D_y)U_k)({\mathcal G}_{j\sigma ks}y)
    =f_{j\sigma\mu}(y) \quad (y\in\gamma_{j\sigma}^\varepsilon).\label{eqBinK}
\end{gather}
Here (and below unless otherwise stated) $j, k=1, \dots, N;$
$\sigma=1, 2;$ $\mu=1, \dots, m;$ $s=0, \dots, S_{j\sigma k}$;
${\mathbf P}_j(y, D_y)$ and $B_{j\sigma\mu ks}(y, D_y)$ are
differential operators of order $2m$ and $m_{j\sigma\mu}$
($m_{j\sigma\mu}\le m-1$), respectively, with $C^\infty$
complex-valued coefficients, i.e.,
$$
 \mathbf P_j(y, D_y)=\sum\limits_{|\alpha|\le 2m}
 p_{j\alpha}(y)D^\alpha_y,\qquad
B_{j\sigma\mu ks}(y, D_y)= \sum\limits_{|\alpha|\le
m_{j\sigma\mu}}
 b_{j\sigma\mu ks\alpha}(y)D^\alpha_y;
$$
${\mathcal G}_{j\sigma ks}$ is the operator of rotation by an
angle~$\omega_{j\sigma ks}$ and of the homothety with a
coefficient~$\chi_{j\sigma ks}$ ($\chi_{j\sigma ks}>0$) in the
$y$-plane. Moreover,
$$
|(-1)^\sigma b_{j}+\omega_{j\sigma
ks}|<b_{k}\qquad\text{for}\qquad (k,s)\ne(j,0)
$$
(cf. Remark~\ref{remK1}) and
$$
\omega_{j\sigma j0}=0,\qquad \chi_{j\sigma j0}=1
$$
(i.e., ${\mathcal G}_{j\sigma j0}y\equiv y$).

%% file: sect2.tex
\section{The Fredholm Property of Nonlocal Problems}\label{sectProof}

In this section, we prove the following result.

\begin{theorem}\label{thP_BFred}
Let the operator ${\bf P}(y, D_y)$ be properly elliptic on
$\overline{G}$, and let the system $\{B_{i\mu 0}(y,
D_y)\}_{\mu=1}^m$ satisfy the Lopatinsky condition on the curve
$\overline{\Gamma_i}$ with respect to ${\bf P}(y, D_y)$ for all
$i=1, \dots, N$. Assume that Conditions~$\ref{condK1}$
and~$\ref{condSeparK23}$ are fulfilled. Then the operator $\mathbf
P$ has the Fredholm property.
\end{theorem}

\begin{remark} One can assign a bounded operator
(acting from $W^{2m}(G)$ to $L_2(G)$) to problem~(\ref{eqPinG}),
(\ref{eqBinG}). Such an operator is studied
in~\cite{GurRJMP03,GurRJMP04}; it is proved that, unlike the case
treated in the present paper, whether or not the bounded operator
has the Fredholm property depends both on spectral properties of
auxiliary nonlocal problems with a parameter and on the validity
of some algebraic relations between the operators ${\bf P}(y,
D_y)$, $\mathbf B_{i\mu}^0$, and $\mathbf B_{i\mu}^1$ at the
points of the set $\mathcal K$.
\end{remark}

\subsection{Finite Dimensionality of the Kernel}

In this subsection, we prove that the kernel of the operator
$\mathbf P$ is of finite dimension. To do this, we preliminarily
study the smoothness of generalized solution of
problem~(\ref{eqPinG}), (\ref{eqBinG}). We first study the
smoothness outside a neighborhood of the set $\mathcal K$ and then
near $\mathcal K$. The following lemma generalizes part~1 of
Theorem~5 in~\cite{SkDu82}.

\begin{lemma}\label{lSmoothOutsideK}
Let Condition~$\ref{condSeparK23}$ hold, and let $u\in W^{m}(G)$
be a generalized solution of problem~\eqref{eqPinG},
\eqref{eqBinG} with right-hand side $f_0\in L_2(G)$. Then
\begin{equation}\label{eqSmoothOutsideK}
u\in W^{2m}\bigl(G\setminus\overline{\mathcal O_\delta(\mathcal
K)}\bigr)\qquad \text{for any}\qquad \delta>0.
\end{equation}
\end{lemma}
\begin{proof}
1) Denote by $W^k_\loc(G)$ the set of distributions $v$ on $G$
such that $\psi v\in W^k(G)$ for all $\psi\in C_0^\infty(G)$. It
follows from Theorem~3.2 in~\cite[Chap.~2]{LM} that
\begin{equation}\label{eqSmoothLocG}
 u\in W^{2m}_\loc(G).
\end{equation}
This relation and estimate~\eqref{eqSeparK23''} imply that
\begin{equation}\label{eqSmoothOutsideK_2}
\mathbf B_{i\mu}^2 u\in
W^{2m-m_{i\mu}-1/2}(\Gamma_i\setminus\overline{\mathcal
O_{\varkappa_2}(\mathcal K)}).
\end{equation}

Fix an arbitrary point $g\in\Gamma_i\setminus\overline{\mathcal
O_{\varkappa_2}(\mathcal K)}$. Take a number $\delta>0$ such that
\begin{equation}\label{eqSmoothOutsideK_3}
 \overline{\mathcal
O_{\delta}(g)\cap\Gamma_i}\subset\Gamma_i\setminus\overline{\mathcal
O_{\varkappa_2}(\mathcal K)}.
\end{equation}

Then the function $u$ is a solution of the following ``local''
problem in the neighborhood $\mathcal O_{\delta}(g)$:
\begin{align}
\mathbf P(y, D_y)u&=f_0(y)\quad (y\in\mathcal
O_{\delta}(g)\cap G),\label{eqSmoothOutsideK_4}\\
B_{i\mu0}(y, D_y)u&=f'_{i\mu}(y)\quad (y\in\mathcal
O_{\delta}(g)\cap \Gamma_i;\ \mu=1, \dots,
m),\label{eqSmoothOutsideK_5}
\end{align}
where $f'_{i\mu}(y)=-\mathbf B_{i\mu}^1u(y)-\mathbf
B_{i\mu}^2u(y)$ for $y\in\mathcal O_{\delta}(g)\cap \Gamma_i$. It
follows from relations~\eqref{eqSmoothLocG},
\eqref{eqSmoothOutsideK_2}, and~\eqref{eqSmoothOutsideK_3} and
from the definition of the operator $\mathbf B_{i\mu}^1$ that
$f'_{i\mu}\in W^{2m-m_{i\mu}-1/2}(\mathcal O_{\delta}(g)\cap
\Gamma_i)$.

Applying Theorem~8.2 in~\cite[Chap.~2]{LM}\footnote{It is
additionally supposed in Theorem~8.2 in~\cite[Chap.~2]{LM} that
the operators $B_{i\mu0}(y, D_y)$ are normal on $\Gamma_i$, while
their orders are not equal to one another. However, it is easy to
check that the theorem mentioned remains valid without these
assumptions (see~\cite[Chap.~2, Sec.~8.3]{LM}).} to
problem~(\ref{eqSmoothOutsideK_4}), (\ref{eqSmoothOutsideK_5}), we
obtain
\begin{equation}\label{eqSmoothOutsideK_6}
 u\in W^{2m}(\mathcal O_{\delta/2}(g)\cap G).
\end{equation}

By using a partition of unity, we infer from~(\ref{eqSmoothLocG})
and~(\ref{eqSmoothOutsideK_6}) that
\begin{equation}\label{eqSmoothOutsideK_7}
u\in W^{2m}\bigl(G\setminus\overline{\mathcal
O_{\varkappa_1}(\mathcal K)}\bigr).
\end{equation}

\smallskip

2) It follows from the belonging~\eqref{eqSmoothOutsideK_7} and
from inequality~\eqref{eqSeparK23'} that
\begin{equation}\label{eqSmoothOutsideK_8}
\mathbf B_{i\mu}^2 u\in W^{2m-m_{i\mu}-1/2}(\Gamma_i).
\end{equation}

Taking into account~\eqref{eqSmoothOutsideK_8}, we can repeat the
arguments of part~1) of this proof for arbitrary $g\in\Gamma_i$
and $\delta>0$ such that
\begin{equation}\notag
 \overline{\mathcal
O_{\delta}(g)\cap\Gamma_i}\subset\Gamma_i.
\end{equation}
As a result, we obtain the belonging~\eqref{eqSmoothOutsideK_6}
valid for an arbitrary point $g\in\Gamma_i$. Combining this fact
with relation~(\ref{eqSmoothLocG}) and using a partition of unity,
we deduce~\eqref{eqSmoothOutsideK}.
\end{proof}

Now we study the smoothness of solutions of
problem~(\ref{eqPinG}), (\ref{eqBinG}) in a neighborhood of the
set $\mathcal K$. Since generalized solutions can have power-law
singularities near the set $\mathcal K$ (see~\cite{SkMs86}), it is
natural to consider these solutions in weighted spaces. Let us
introduce these spaces.

Assume that either $Q=\{y\in{\mathbb R}^2:\ r>0,\ |\omega|<b\}$ or
$Q=\{y\in{\mathbb R}^2:\ 0<r<d,\ |\omega|<b\}$, $0<b<\pi$, $d>0$,
or $Q=G$. In the first and second cases, we set $\mathcal
M=\{0\}$, while in the third case we set $\mathcal M=\mathcal K$.
Introduce the space $H_a^k(Q)$ as the completion of the set
$C_0^\infty(\overline{ Q}\setminus \mathcal M)$ with respect to
the norm
$$
 \|w\|_{H_a^k(Q)}=\left(
    \sum_{|\alpha|\le k}\int\limits_Q \rho^{2(a-k+|\alpha|)} |D_y^\alpha w|^2 dy
                                       \right)^{1/2},
$$
where $a\in \mathbb R$, $k\ge 0$ is an integer, and
$\rho=\rho(y)=\dist(y,\mathcal M)$. For integer $k\ge1$, denote by
$H_a^{k-1/2}(\gamma)$ the space of traces on a smooth curve
$\gamma\subset\overline{ Q}$ with the norm
\begin{equation}\label{eqTraceNormH}
\|\psi\|_{H_a^{k-1/2}(\gamma)}=\inf\|w\|_{H_a^k(Q)} \quad (w\in
H_a^k(Q): w|_\gamma = \psi).
\end{equation}

\medskip

Let $u$ be a generalized solution of problem~\eqref{eqPinG},
\eqref{eqBinG}, and let $U_j(y')=u_j(y(y'))$, $j=1,\dots,N$, be
the functions corresponding to the set (orbit) $\mathcal K$ and
satisfying problem~\eqref{eqPinK}, \eqref{eqBinK} with right-hand
side $\{f_j, f_{j\sigma\mu}\}$ (see
Sec.~\ref{subsectStatementNearK}).

Set
$$
d_1=\min\{\chi_{j\sigma ks}, 1\}/2,\qquad d_2=2\max\{\chi_{j\sigma
ks}, 1\}.
$$
Take a sufficiently small $\varepsilon$ such that
$d_2\varepsilon<\varepsilon_1$. It follows from
Lemma~\ref{lSmoothOutsideK} that
\begin{equation}\label{eqP_BFiniteKer5'}
U_j\in
W^{2m}(K_{j}^{d_2\varepsilon}\cap\{|y|>\delta\})\qquad\text{for
any}\qquad \delta>0.
\end{equation}
Further, it follows from the belonging $U_j\in
W^m(K_{j}^{d_2\varepsilon})$ and from Lemma~5.2 in~\cite{KovSk}
that
\begin{equation}\label{eqP_BFiniteKer5}
U_j\in H_{a-m}^m(K_{j}^{d_2\varepsilon})\subset
H_{a-2m}^0(K_{j}^{d_2\varepsilon}),\quad a>2m-1.
\end{equation}
Finally,  $f_j\in L_2(K_{j}^{\varepsilon})$ and, by virtue of
Lemma~\ref{lSmoothOutsideK} and estimate~\eqref{eqSeparK23'},
$f_{j\sigma\mu}\in
W^{2m-m_{j\sigma\mu}-1/2}(\gamma_{j\sigma}^\varepsilon)$.
Therefore, by Lemma~5.2 in~\cite{KovSk},
\begin{equation}\label{eqP_BFiniteKer4}
f_j\in H_a^{0}(K_{j}^{\varepsilon}),\quad f_{j\sigma\mu}\in
H_a^{2m-m_{j\sigma\mu}-1/2}(\gamma_{j\sigma}^\varepsilon),\quad
a>2m-1.
\end{equation}

The following two lemmas enable us to prove that $U_j\in
H_{a}^{2m}(K_{j}^{\varepsilon/d_2^3})$ whenever
relations~(\ref{eqP_BFiniteKer5'})--(\ref{eqP_BFiniteKer4}) hold.

\smallskip

Set
$$
K_{jq}=K_j\cap\{\varepsilon d_2^{-3}d_1^{4-q}/2<|y|<\varepsilon
d_2^{-3}d_2^{4-q}\},\qquad q=0, \dots, 4.
$$
\begin{lemma}\label{lAprinKjq}
Let Condition~$\ref{condK1}$ hold. Then the estimate
 \begin{multline}\label{eqAprinKjq}
  \sum\limits_j\|U_j\|_{W^{2m}(K_{j4})}\le c\sum\limits_j\bigl\{\|\mathbf
  P_j(y, D_y)U_j\|_{L_2(K_{j1})}\\
  +\sum\limits_{\sigma,\mu}\|\mathbf
  B_{j\sigma\mu}(y, D_y)U|_{\gamma_{j\sigma}\cap\overline{ K_{j1}}}
  \|_{W^{2m-m_{j\sigma\mu}-1/2}(\gamma_{j\sigma}\cap\overline{ K_{j1}})}
    +\|U_j\|_{L_2(K_{j1})}\big\}
 \end{multline}
holds for any $U\in\prod\limits_j W^{2m}(K_{j0})$, where $c>0$
does not depend on $U$.
\end{lemma}
\begin{proof}
It follows from the general theory of elliptic problems that
\begin{multline}\label{eqAprinKjq1}
\|U_j\|_{W^{2m}(K_{j4})} \le k_1\big(\|\mathbf
  P_j(y, D_y)U_j\|_{L_2(K_{j3})}\\
  +\sum\limits_{\sigma,\mu}\|
  B_{j\sigma\mu j0}(y, D_y)U_j|_{\gamma_{j\sigma}\cap\overline{ K_{j3}}}
  \|_{W^{2m-m_{j\sigma\mu}-1/2}(\gamma_{j\sigma}\cap\overline{ K_{j3}})}
    +\|U_j\|_{L_2(K_{j3})}).
\end{multline}
Let $(k,s)\ne(j,0)$; then the set ${\mathcal G}_{j\sigma
ks}(\gamma_{j\sigma})\cap\overline{ K_{k2}}$ lies strictly inside
the domain $K_{k1}$. Therefore, using the boundedness of the trace
operator on the corresponding Sobolev spaces, we obtain (similarly
to~\eqref{eqAprinKjq1})
\begin{multline}\label{eqAprinKjq2}
\|B_{j\sigma\mu ks}(y, D_y)U_k({\mathcal G}_{j\sigma
ks}y)|_{\gamma_{j\sigma}\cap\overline{
K_{j3}}}\|_{W^{2m-m_{j\sigma\mu}-1/2}(\gamma_{j\sigma}\cap\overline{
K_{j3}})}\\ \le k_2\|B_{j\sigma\mu ks}(y, D_y)U_k|_{{\mathcal
G}_{j\sigma ks}(\gamma_{j\sigma})\cap\overline{
K_{k2}}}\|_{W^{2m-m_{j\sigma\mu}-1/2}({\mathcal G}_{j\sigma
ks}(\gamma_{j\sigma})\cap\overline{ K_{k2}})}\\
\le k_3 (\|\mathbf
  P_j(y, D_y)U_k\|_{L_2(K_{k1})}+\|U_k\|_{L_2(K_{k1})}).
\end{multline}
Estimates~\eqref{eqAprinKjq1} and~\eqref{eqAprinKjq2}
imply~\eqref{eqAprinKjq}.
\end{proof}

\begin{remark}\label{remConstinApr}
Assume that the norm (in $C^0(\overline{ K_{j1}})$) of the
coefficients $p_{j\alpha}$ of the operators $\mathbf P_j(y, D_y)$
and the norms (in $C^{2m-m_{j\sigma\mu}}(\overline{ K_{j0}})$) of
the coefficients $b_{j\sigma\mu ks\alpha}$ of the operators
$B_{j\sigma\mu ks}(y, D_y)$ do not exceed some constant $C$. Let
the norms (in $C^1(\overline{ K_{j1}})$) of the coefficients
$p_{j\alpha}$, $|\alpha|=2m$, at senior terms of the operators
$\mathbf P_j(y, D_y)$ not exceed the same constant $C$. In that
case, the constant $c$ occurring in inequality~\eqref{eqAprinKjq}
depends only on $C$, on the constant $A$ in~\eqref{eqPEllipinG},
and on the constant $D$ in~\eqref{eqCover}.
\end{remark}

\begin{lemma}\label{lAprinK}
Let Condition~$\ref{condK1}$ hold. Assume that a function $U$
satisfies relations~\eqref{eqP_BFiniteKer5'}
and~\eqref{eqP_BFiniteKer5} and is a solution of
problem~\eqref{eqPinK}, \eqref{eqBinK} with right-hand side $
\{f_j, f_{j\sigma\mu}\}$ satisfying
relations~\eqref{eqP_BFiniteKer4}. Then $U\in \prod\limits_j
H_a^{2m}(K_{j}^{\varepsilon/d_2^3})$ and
\begin{equation}\label{eqAprinK}
\sum\limits_j \|U_j\|_{H_a^{2m}(K_{j}^{\varepsilon/d_2^3})}\le
c\sum\limits_j\big\{\|f_j\|_{H_a^{0}(K_{j}^{\varepsilon})}+
\sum\limits_{\sigma,\mu}\|f_{j\sigma\mu}\|_{H_a^{2m-m_{j\sigma\mu}-1/2}
(\gamma_{j\sigma}^\varepsilon)}
+\|U_j\|_{H_{a-2m}^{0}(K_{j}^{\varepsilon})}\big\},
\end{equation}
where $c>0$ does not depend on $U$.
\end{lemma}
\begin{proof}
Set
$$
K_{jq}^s=K_j\cap\{\varepsilon d_2^{-3}
d_1^{4-q}2^{-s-1}<|y|<\varepsilon d_2^{-3}d_2^{4-q}2^{-s}\},\qquad
s=0,1,2, \dots\,.
$$
Clearly,
\begin{equation}\label{eqKjEpsilon}
\bigcup_{s=0}^{\infty}K_{j1}^s=K_{j}^{\varepsilon},\qquad
\bigcup_{s=0}^{\infty}K_{j4}^s=K_{j}^{\varepsilon/d_2^3}.
\end{equation}

Set $U_j^s(y')=U_j(2^{-s}y')$ and make the change of variables
$y=2^{-s}y'$ in the equation
$$
 \mathbf P_j(y, D_y)U_j\equiv\sum\limits_{|\alpha|\le 2m}
 p_{j\alpha}(y)D^\alpha_y U_j(y)=f_j(y)\quad (y\in K_{j1}^s)
$$
and in the nonlocal conditions
$$
\sum\limits_{k,s}\sum\limits_{|\alpha|\le m_{j\sigma\mu}}
 b_{j\sigma\mu ks\alpha}(x)D^\alpha_x U_j(x)|_{x={\mathcal G}_{j\sigma
 ks}y}=f_{j\sigma\mu}(y)\quad (y\in \gamma_{j\sigma}\cap\overline{
 K_{j1}^s});
$$
multiplying the first equation obtained by $2^{-s\cdot 2m}$ and
the second one by $2^{-s\cdot m_{j\sigma\mu}}$, we have
\begin{equation}\label{eqAprinK1}
 \sum\limits_{|\alpha|\le 2m}
 p_{j\alpha}^s(y')2^{s(|\alpha|-2m)}D^\alpha_{y'}
 U^s_j(y')=2^{-s\cdot2m}f^s_j(y')\quad (y'\in K_{j1}^0),
\end{equation}
\begin{equation}\label{eqAprinK1'}
\sum\limits_{k,s}\sum\limits_{|\alpha|\le m_{j\sigma\mu}}
 b^s_{j\sigma\mu ks\alpha}(x')2^{s(|\alpha|-m_{j\sigma\mu})}D^\alpha_{x'} U^s_j(x')|_{x'={\mathcal G}_{j\sigma
 ks}y'}
=2^{-s\cdot m_{j\sigma\mu}}f_{j\sigma\mu}^s(y')\quad (y'\in
\gamma_{j\sigma}\cap\overline{
 K_{j1}^0}),
\end{equation}
where
$$
p_{j\alpha}^s(y')=p_{j\alpha}(2^{-s}y'),\qquad b^s_{j\sigma\mu ks
\alpha}(x')=b_{j\sigma\mu ks\alpha}(2^{-s}x'),
$$
$$
f_{j}^s(y')=f_{j}(2^{-s}y'),\qquad
f_{j\sigma\mu}^s(y')=f_{j\sigma\mu}(2^{-s}y').
$$

Applying Lemma~\ref{lAprinKjq} to problem~(\ref{eqAprinK1}),
(\ref{eqAprinK1'}), we obtain
\begin{multline}\label{eqAprinK2}
\sum\limits_j\|U_j^s\|_{W^{2m}(K_{j4}^0)}\le k_1\sum\limits_j
\bigl\{\|2^{-s\cdot2m}f_j^s\|_{L_2(K_{j1}^0)}\\
+\sum\limits_{\sigma,\mu} \|2^{-s\cdot
m_{j\sigma\mu}}f_{j\sigma\mu}^s\|_{W^{2m-m_{j\sigma\mu}-1/2}
(\gamma_{j\sigma}\cap \overline{ K_{j1}^0})}
+\|U_j^s\|_{L_2(K_{j1}^0)}\bigr\},
\end{multline}
where $k_1>0$ does not depend on $s$ due to
Remark~\ref{remConstinApr}.

Consider a function $\Phi_{j\sigma\mu}\in
H_a^{2m-m_{j\sigma\mu}}(K_j)$ satisfying the following conditions:
$\Phi_{j\sigma\mu}|_{\gamma_{j\sigma}^\varepsilon}=f_{j\sigma\mu}$
and
\begin{equation}\label{eqAprinK3}
\|\Phi_{j\sigma\mu}\|_{H_a^{2m-m_{j\sigma\mu}}(K_j^\varepsilon)}\le
2\|f_{j\sigma\mu}\|_{
H_a^{2m-m_{j\sigma\mu}-1/2}(\gamma_{j\sigma}^\varepsilon)}
\end{equation}
(the existence of such a function follows
from~\eqref{eqTraceNormH}). Then
$\Phi^s_{j\sigma\mu}|_{\gamma_{j\sigma}\cap
\overline{K_{j1}^0}}=f^s_{j\sigma\mu}$, where
$\Phi^s_{j\sigma\mu}(y')=\Phi_{j\sigma\mu}(2^{-s}y')$. Therefore,
relations~\eqref{eqAprinK2} and~\eqref{eqTraceNormW} imply
\begin{equation}\label{eqAprinK4}
\sum\limits_j\|U_j^s\|_{W^{2m}(K_{j4}^0)}\le k_1\sum\limits_j
\bigl\{\|2^{-s\cdot2m}f_j^s\|_{L_2(K_{j1}^0)}+\sum\limits_{\sigma,\mu}
\|2^{-s\cdot
m_{j\sigma\mu}}\Phi_{j\sigma\mu}^s\|_{W^{2m-m_{j\sigma\mu}}
(K_{j1}^0)} +\|U_j^s\|_{L_2(K_{j1}^0)}\bigr\}.
\end{equation}

Making the inverse change of variables $y'=2^s y$ in
inequality~\eqref{eqAprinK4}, we obtain
\begin{multline}\label{eqAprinK5}
  \sum\limits_j\sum\limits_{|\alpha|\le 2m}
  \|2^{-s|\alpha|}D^\alpha_y U_j\|_{L_2(K_{j4}^s)}\le k_1\sum\limits_j
\bigl\{\|2^{-s\cdot2m}f_j\|_{L_2(K_{j1}^s)}\\+\sum\limits_{\sigma,\mu}
\sum\limits_{|\alpha|\le 2m-m_{j\sigma\mu}}
\|2^{-s(|\alpha|+m_{j\sigma\mu})}\Phi_{j\sigma\mu}\|_{L_2
(K_{j1}^s)} +\|U_j\|_{L_2(K_{j1}^s)}\bigr\}.
\end{multline}
Multiplying inequality~(\ref{eqAprinK5}) by $2^{-s(a-2m)}$,
summing with respect to $s$, and taking into
account~\eqref{eqAprinK3} and~\eqref{eqKjEpsilon}, we
deduce~(\ref{eqAprinK}).
\end{proof}

Combining Lemma~\ref{lAprinK} with Lemma~\ref{lSmoothOutsideK}
yields $u\in H_a^{2m}(G)$, $a>2m-1$, where $u$ is an arbitrary
generalized solution of problem~\eqref{eqPinG}, \eqref{eqBinG}
with the right-hand side $f_0\in L_2(G)$.

It follows from Lemma~2.1 in~\cite{SkDu90} and from Theorem~3.2
in~\cite{SkDu91} that the set of solutions from $H_a^{2m}(G)$ of
problem~\eqref{eqPinG}, \eqref{eqBinG} with right-hand side
$f_0=0$ is of finite dimension for almost all $a>2m-1$. Thus, we
have proved the following result.
\begin{lemma}\label{lP_BFiniteKer}
Let Conditions~$\ref{condK1}$ and~$\ref{condSeparK23}$ hold. Then
the kernel of the operator $\mathbf P$ is of finite dimension.
\end{lemma}

\bigskip

\subsection{Closedness of the Operator and its Image.
Finite Dimensionality of the Cokernel}

To prove that the operator $\mathbf P$ has the Fredholm property,
we need to consider problem~(\ref{eqPinG}), (\ref{eqBinG}) on
weighted spaces with weight $a$ such that $0<a\le m$. Now the
difficulty is that the belonging $u\in H_a^{2m}(G)$ does not imply
that $\mathbf B^2_{i\mu}u\in H_a^{2m-m_{i\mu}-1/2}(\Gamma_i)$;
therefore, the sum
$$
\mathbf B_{i\mu}u=\mathbf B_{i\mu}^0
u+\mathbf B_{i\mu}^1 u+\mathbf B_{i\mu}^2 u
$$
does not necessarily belong to $H_a^{2m-m_{i\mu}-1/2}(\Gamma_i)$.
One can only guarantee  that $\mathbf B_{i\mu}u\in
H_{a'}^{2m-m_{i\mu}-1/2}(\Gamma_i)$, where $a'>2m-1$ (which
follows from the fact that $\mathbf B_{i\mu}u\in
W^{2m-m_{i\mu}-1/2}(\Gamma_i)$ and from Lemma~5.2
in~\cite{KovSk}). However, it is proved
in~\cite[Sec.~6]{GurRJMP04} that
$$
\{{\bf P}(y, D_y)u,\ \mathbf B_{i\mu}u\}\in \mathcal
 H_a^0(G,\Gamma)\dotplus \mathcal R_a^0(G,\Gamma)\qquad\text{for all}\qquad
  u\in H_a^{2m}(G),\ a>0,
$$
where $\mathcal H_a^0(G,\Gamma)=H_a^0(G)\times\prod\limits_{i,\mu}
H_{a}^{2m-m_{i\mu}-1/2}(\Gamma_i)$ and $\mathcal R_a^0(G,\Gamma)$
is some \textit{finite-dimensional} space naturally embedded in
$\{0\}\times\prod\limits_{i,\mu}
H_{a'}^{2m-m_{i\mu}-1/2}(\Gamma_i)$ for any $a'>2m-1$. In
particular, this means that the space $\mathcal R_a^0(G,\Gamma)$
contains only functions of the form $\{0,f_{i\mu}\}$, where
$f_{i\mu}\in H_{a'}^{2m-m_{i\mu}-1/2}(\Gamma_i)$ and
$f_{i\mu}\notin H_{a}^{2m-m_{i\mu}-1/2}(\Gamma_i)$. Fix some
$a'>2m-1$. Then any function
$$
\{f_0,f_{i\mu}\}\in\mathcal
 H_a^0(G,\Gamma)\dotplus \mathcal R_a^0(G,\Gamma)
$$
can be represented as follows:
$$\{f_0,f_{i\mu}\}=\{f_0,f_{i\mu}^1\}+\{0,f_{i\mu}^2\},$$ where
$\{f_0,f_{i\mu}^1\}\in\mathcal H_a^0(G,\Gamma)$ and
$\{0,f_{i\mu}^2\}\in\mathcal R_a^0(G,\Gamma)$, and its norm is
given by
$$
\|\{f_0,f_{i\mu}\}\|_{\mathcal H_a^0(G,\Gamma)\dotplus \mathcal
R_a^0(G,\Gamma)}=\Big(\|\{f_0,f_{i\mu}^1\}\|^2_{\mathcal
 H_a^0(G,\Gamma)}+\sum\limits_{i,\mu}\|f_{i\mu}^2\|^2_{H_{a'}^{2m-m_{i\mu}-1/2}(\Gamma_i)}\Big)^{1/2}.
$$

Furthermore, it follows from Theorem~6.1 in~\cite{GurRJMP04} that
the operator
\begin{equation}\notag
 \mathbf L_a=\{{\bf P}(y, D_y),\ \mathbf B_{i\mu}\}: H_a^{2m}(G)\to  \mathcal
 H_a^0(G,\Gamma)\dotplus\mathcal R_a^0(G,\Gamma),\quad a>0,
\end{equation}
has the Fredholm property for almost all $a>0$. In other words, if
$u\in\ H_a^{2m}(G)$, then $\mathbf L_au$ ``belongs'' to the space
$\mathcal H_a^0(G,\Gamma)$ up to a function of the form
$\{0,f_{i\mu}\}$ from the {\it finite-dimensional} space $\mathcal
R_a^0(G,\Gamma)$.

Using the Fredholm property of the operator $\mathbf L_a$, we
prove the following result.

\begin{lemma}\label{lP_BClosed}
Let Conditions~$\ref{condK1}$ and~$\ref{condSeparK23}$ hold. Then
the operator $\mathbf P$ is closed, its image $\mathcal R(\mathbf
P)$ is closed, and $\codim\mathcal R(\mathbf P)<\infty$.
\end{lemma}
\begin{proof}
1) Let $0<a\le m$. We consider the auxiliary unbounded operator
$$\mathbf P_{a}:\Dom(\mathbf P_{a})\subset L_2(G)\to L_2(G)$$ given
by
$$
\mathbf P_{a}u=\mathbf P(y, D_y)u,\quad u\in \Dom(\mathbf
P_{a})=\{u\in H_a^{2m}(G):\ \mathbf B_{i\mu}u=0,\ \mathbf P(y,
D_y)u\in L_2(G)\}.
$$
Fix a number $a$, $0<a\le m$, such that the operator ${\mathbf
L}_a$ has the Fredholm property. Let us show that the operator
$\mathbf P_{a}$ also has the Fredholm property.

Since ${\mathbf L}_a$ has the Fredholm property, it follows from
the compactness of the embedding $H_a^{2m}(G)\subset H_a^0(G)$
(see Lemma~3.5 in~\cite{KondrTMMO67}) and from Theorem~7.1
in~\cite{Kr} that
\begin{equation}\label{eqP_BClosed1'}
 \|u\|_{H_a^{2m}(G)}\le k_1(\|{\mathbf L}_au\|_{\mathcal
H_a^0(G,\Gamma)\dotplus\mathcal R_a^0(G,\Gamma)}+
 \|u\|_{H_a^{0}(G)})
\end{equation}
for all  $u\in H_a^{2m}(G)$.

Now we take a function $u\in \Dom(\mathbf P_{a})$. Then ${\mathbf
L}_au=\{\mathbf P(y, D_y)u,0\}$, $\mathbf P(y, D_y)u\in
L_2(G)\subset H_a^0(G)$, and hence
$$
\|{\mathbf L}_au\|_{\mathcal H_a^0(G,\Gamma)\dotplus\mathcal
R_a^0(G,\Gamma)}=\|\mathbf P(y, D_y)u\|_{H_a^0(G)}.
$$
Combining this relation with~\eqref{eqP_BClosed1'} and taking into
account the boundedness of the embedding $L_2(G)\subset H_a^0(G)$
for $a>0$, we obtain
\begin{equation}\label{eqP_BClosed1}
 \|u\|_{H_a^{2m}(G)}\le k_2(\|\mathbf P(y, D_y)u\|_{H_a^0(G)}+
 \|u\|_{H_a^{0}(G)})
 \le k_3(\|\mathbf P(y, D_y)u\|_{L_2(G)}+
 \|u\|_{L_2(G)}),
\end{equation}
where $u\in \Dom(\mathbf P_{a})$. It follows from
inequality~(\ref{eqP_BClosed1}) that the operator $\mathbf P_{a}$
is closed. Therefore, using~(\ref{eqP_BClosed1}) and applying
Lemma~7.1 in~\cite{Kr} again, we obtain that $ \dim\ker\mathbf
P_{a}<\infty $ (clearly, $\ker\mathbf P_{a}=\ker\mathbf L_{a}$)
and the image $\mathcal R(\mathbf P_{a})$ is closed.

Consider an arbitrary function $f_0\in L_2(G)$. Clearly, $f_0\in
H_a^0(G)$. By Corollary~6.1 in~\cite{GurRJMP04}, there exist
functionals $F_1, \dots, F_{q_0}$ from the adjoint space
${\mathcal H}_a^0(G,\Gamma)^*$ such that problem~(\ref{eqPinG}),
(\ref{eqBinG}) admits a solution $u\in H_a^{2m}(G)$ whenever
$$
\langle\{f_0,0\},\ F_q\rangle=0,\qquad q=1, \dots, q_0.
$$
Since
$$
 |\langle\{f_0,0\},\ F_q\rangle|\le k_4\|f_0\|_{H_a^0(G)}\le k_5\|f_0\|_{L_2(G)},
$$
it follows from Riesz' theorem on the general form of a continuous
linear functional on a Hilbert space that there exist functions
$f_1, \dots, f_{q_0}\in L_2(G)$ such that
$$
\langle\{f_0,0\},\ F_q\rangle=(f_0,\ f_q)_{L_2(G)},\qquad q=1,
\dots, q_0.
$$
Therefore, $\codim\mathcal R(\mathbf P_{a})\le q_0$.

Thus, we have proved that the operator $\mathbf P_{a}$ has the
Fredholm property.

2) Since $H_a^{2m}(G)\subset H_{a-m}^m(G)\subset W^m(G)$ for $a\le
m$, it follows that
\begin{equation}\label{eqP_BClosed2}
 \mathbf P_{a}\subset \mathbf P.
\end{equation}
In particular, relation~(\ref{eqP_BClosed2}) implies that the
image $\mathcal R(\mathbf P)$ is closed and $$\codim\mathcal
R(\mathbf P)\le\codim\mathcal R(\mathbf P_{a})\le q_0.$$

It remains to prove that the operator $\mathbf P$ is
closed.\footnote{Note that the closedness of the image of some
operator $\mathbf P$ on a Hilbert space and the finite
dimensionality of its kernel and cokernel do not imply the
closedness of $\mathbf P$ itself; this can be shown by using
arguments close to that in~\cite[Chap.~2, Sec.~18]{AhGl}. However,
if we additionally suppose that the operator $\mathbf P$ is an
extension of a Fredholm operator, then we prove that $\mathbf P$
is closed.} Denote by $h_1,\dots,h_k$ some basis of the space
$$
\mathcal R(\mathbf P_{a})^\bot=\mathcal R(\mathbf
P)\ominus\mathcal R(\mathbf P_{a}).
$$
Then there exist functions $v_1,\dots,v_k\in\Dom(\mathbf P)$ such
that $\mathbf Pv_j=h_j$, $j=1,\dots,k$. Since $h_j\notin\mathcal
R(\mathbf P_{a})$, it follows that $v_j\notin\Dom(\mathbf P_{a})$.
It is also clear that the functions $v_1,\dots,v_k$ are linearly
independent because the functions $h_1,\dots,h_k$ have this
property.

Consider the finite-dimensional space
$$
\mathcal N=\Span(v_1,\dots,v_k,\ker\mathbf P)\ominus\ker\mathbf
P_{a}.
$$
It is easy to see that $\mathcal N\cap\Dom\mathbf P_{a}=\{0\}$.
Indeed, if $u\in\mathcal N\cap\Dom\mathbf P_{a}$, then
$$
u=\sum\limits_{i=1}^k\alpha_iv_i+v,
$$
where $\alpha_i$ are some constants and $v\in\ker\mathbf P$.
Therefore, taking into account~\eqref{eqP_BClosed2}, we have
$$
\sum\limits_{i=1}^k\alpha_ih_i=\mathbf Pu=\mathbf P_au\in\mathcal
R(\mathbf P_{a}).
$$
Hence, $\alpha_i=0$, $i=1,\dots,k$, which implies that $u=v$.
Using~\eqref{eqP_BClosed2} again, we see that $u=v\in\ker\mathbf
P_{a}$. Combining this fact with the definition of the space
$\mathcal N$ yields $u=0$.

Let $\Gr\mathbf P$ ($\Gr\mathbf P_{a}$) denote the graph of the
operator $\mathbf P$ ($\mathbf P_{a}$). As is known, the operator
$\mathbf P$ ($\mathbf P_{a}$) is closed if and only if its graph
$\Gr\mathbf P$ ($\Gr\mathbf P_{a}$) is closed in $L_2(G)\times
L_2(G)$.

Note that $\Gr\mathbf P_{a}$ is closed (as the graph of the closed
operator) and $\Gr\mathbf P_{a}\subset\Gr\mathbf P$, while the
spaces $\mathcal N$ and $\mathcal R(\mathbf P_{a})^\bot$ are of
finite dimension. Therefore, to prove that the operator $\mathbf
P$ is closed, it suffices to show that
\begin{equation}\label{eqP_BClosed3}
\Gr\mathbf P\subset\Gr\mathbf P_{a}\dotplus(\mathcal N\times
\mathcal R(\mathbf P_{a})^\bot).
\end{equation}
Clearly, the sum in~\eqref{eqP_BClosed3} is direct. Indeed, if
$$(u,f)\in\Gr\mathbf P_{a}\cap(\mathcal N\times \mathcal R(\mathbf
P_{a})^\bot),$$ then $u\in\Dom\mathbf P_{a}\cap\mathcal N=\{0\}$,
and hence $(u,f)=(u,\mathbf P_{a}u)=(0,0)$.

Further, let $(u,f)\in\Gr\mathbf P$, i.e., $u\in\Dom\mathbf P$ and
$f=\mathbf Pu$. We represent the function $f$ as follows:
$$
f=f_1+ f_2,
$$
where $f_1\in\mathcal R(\mathbf P_{a})$ and $f_2\in\mathcal
R(\mathbf P_{a})^\bot$. Take an element $u_1\in\Dom(\mathbf
P_{a})$ such that $\mathbf P_{a}u_1=f_1$. Then
$u_2=u-u_1\in\Dom(\mathbf P)$ and $\mathbf Pu_2=f_2$. Without loss
of generality, one can assume that
\begin{equation}\label{eqP_BClosed4}
u_2\,\bot\,\ker\mathbf P_{a};
\end{equation}
if this relation fails, one must take the projection $u_{2a}$ of
the element $u_2$ to $\ker\mathbf P_{a}$ and replace $u_1$ by
$u_1+u_{2a}$ and $u_2$ by $u_2-u_{2a}$. Clearly,
$(u_1,f_1)\in\Gr\mathbf P_{a}$ and, due to~\eqref{eqP_BClosed4},
$(u_2,f_2)\in\mathcal N\times\mathcal R(\mathbf P_{a})^\bot$.

Thus, we have proved relation~\eqref{eqP_BClosed3}, and the lemma
is true.
\end{proof}
Lemmas~\ref{lP_BFiniteKer} and~\ref{lP_BClosed} imply
Theorem~\ref{thP_BFred}.

\begin{remark}
Using results in~\cite{GurIzvRAN}, one can prove that
Theorem~\ref{thP_BFred} remains valid if the transformations
$\Omega_{is}$ are {\it nonlinear} near the points of the set
$\mathcal K$, while the linear parts of $\Omega_{is}$ satisfy
Condition~\ref{condK1} at the points of $\mathcal K$.
\end{remark}
\medskip

The author is grateful to Professor A.~L.~Skubachevskii for
attention to this work.